\documentstyle[11pt]{article}
\newcommand{\beq}{\begin{equation}}
\newcommand{\eeq}{\end{equation}}
\newcommand{\W}{\hphantom{0}}
\newcommand{\fl}{{\rm{fl}}}
\newcommand{\etal}{{\em{et~al}}}

\begin{document}

\title{On the Precision Attainable with Various\\
Floating-Point Number Systems%
\thanks{Earlier versions 
appeared as Report TR RC 3751, IBM Research (February 1972); and in
{\em IEEE Transactions on Computers}, {\bf{C-22}} (June 1973), 601--607
(manuscript received May 15, 1972; revised May 31, 1972).~~Retyped 
(with corrections) in \LaTeX\ by 
Frances Page, %
July 2000.
}}

\author{Richard P.~Brent%
\thanks{\hbox{Former address: Mathematical Sciences Department, IBM
T.J.~Watson Research Center, Yorktown Heights,}
\hbox{N.Y.~10598.~~Current address: Mathematical Sciences Institute,
Australian National University, Canberra.}
\hbox{Copyright \copyright\ 1972--2010, Richard P.~Brent.}
\hspace*{\fill} \hbox{rpb017 typeset using \LaTeX.}}
}

\date{}

\maketitle

\begin{abstract} \normalsize \noindent
For scientific computations on a digital computer the set of real
numbers is usually approximated by a finite set $F$ of ``floating-point''
numbers.  We compare the numerical accuracy possible with different
choices of $F$ having approximately the same range and requiring the
same word length.  In particular, we compare different choices of base
(or radix) in the usual floating-point systems.  The emphasis is on
the choice of $F$, not on the details of the number representation
or the arithmetic, but both rounded and truncated arithmetic are
considered. Theoretical results are given, and some simulations of
typical floating point-computations (forming sums, solving systems of
linear equations, finding eigenvalues) are described.  If the leading
fraction bit of a normalized base~$2$ number is not stored explicitly
(saving a bit), and the criterion is to minimise the mean square
roundoff error, then base~$2$ is best.  If unnormalized numbers are
allowed, so the first bit must be stored explicitly, then base~$4$
(or sometimes base~$8$) is the best of the usual systems.
\end{abstract}

\begin{quote}
\end{quote}

\thispagestyle{empty}

\section{Introduction}

A real number $x$ is usually approximated in a digital computer by an
element ${\fl}(x)$ of a finite set $F$ of ``floating-point'' numbers.
We regard the elements of $F$ as exactly representable real numbers,
and take ${\fl}(x)$ as the floating-point number closest to $x$.  The
definition of ``closest'', rules for breaking ties, and the possibility
of truncating instead of rounding are discussed later.

We restrict our attention to binary computers in which floating-point
numbers are represented in a word (or multiple word) of fixed length
$w$ bits, using some convenient (possibly redundant) code.  Usually
$F$ is a set of numbers of the form
\beq
s \sum^{t}_{i=1} d_i \beta^{e-i}
\eeq

\noindent
where $\beta = 2^k > 1$ is the base (or radix), $t > 0$ is the number
of digits, $s = \pm 1$ is a sign, $e$ is an exponent in some fixed 
range
\beq
m \ < e  \ \leq \ M\,,
\eeq

\noindent
and each $d_i$ is a $\beta$-ary digit $0, 1, \dots, \beta - 1\,$.
Other possible floating-point number systems (i.e, choices of $F$) 
are mentioned in Section~3.

Since the coding of the exponent $e$ and the signed fraction
$(s; d_1, \dots, d_t)$ must fit into $w$ bits, there is a tradeoff
between precision and range.  (A discussion of precision and range
requirements for general scientific computing may be found in
Cody~\cite{7,8}.)  We do not consider this tradeoff; instead we suppose
that the range and word length is prescribed, and we study the
dependence of the precision on the base $\beta$.

With higher bases less bits are needed for the exponent, so more are
available for the fraction (see Section~2 for details).  However, 
more leading fraction bits may be zero, so the best choice of base
is not immediately obvious.  Our aim is to compare the attainable 
precision of systems with different bases.  Theoretical results are
given in Sections 4 and 5, and some simulations are described in
Sections 6 and 7.  The conclusions are summarised in Section~7.

Since we are interested in the precision attainable with different 
number systems, we assume that the arithmetic is the best possible.
In other words, if $x, y \in F$, and $\dagger$ is an arithmetic
operation, we asume that $x \dagger y$ is found to sufficient
accuracy to give the correct (rounded) result ${\fl}(x \dagger y)$.
Ensuring this may be too expensive in practice, but our conclusions
should be valid provided several guard units are used when computing
${\fl}(x \dagger y)$.  The reduction in precision caused by using only a 
small number of guard digits is discussed by Kuki and Cody~\cite{18}.

\section{The Usual Systems}

\setcounter{equation}{0}

A floating-point number of the form (1.1) may be written as
\beq
s \sum^{u}_{j=1} b_j 2^{ke-j}
\eeq

\noindent
where $b_{k(i-1)+1} \cdots b_{ki}$ is the binary form of the $2^k$-ary
digit $d_i$, and $u=kt$ is the number of bits required to code
$d_i, \dots, d_t$.  We use (2.1) in preference to (1.1), and do 
not insist that $t$ must be an integer.  The details of the coding of
the exponent $e$ and the signed fraction $(s; b_1, \dots, b_u)$
in a $w$-bit word do not concern us.

The representation (2.1) is said to be {\it normalised\/} if at least
one of $b_1, \dots, b_k$ is nonzero.  From (2.1) and the bound (1.2)
on $e$, the largest and smallest floating-point numbers having a
normalised representation are
\beq
f_{\max} \ = \ 2^{kM}(1 - 2^{-u})
\eeq
\noindent and
\beq
\hspace*{-13mm}f_{\min} \ =  \ 2^{km} \,,
\eeq

\noindent
respectively.  If the {\it range\/} $R$ of the system is defined to be
$\log_2 (f_{\max}/f_{\min})$ then, negelecting the term $2^{-u}$ in (2.2)
\beq
\hspace*{-25mm}k(M - m) \ = \ R \,.
\eeq

\noindent
Thus, for systems with the same range, $k(M-m)$ is invatiant.

Goldberg~\cite{10}, McKeeman~\cite{21}, and others have observed 
that with base~$2$ the leading fraction bit $b_1$ can be implicit,
provided only normalized representations of nonzero numbers are
allowed and a special exponent is reserved for zero.  Define
\beq
p \ = \ \left\{ \begin{array}{ll}
            2,&\hbox{\hspace*{4mm}if this ``implicit-first-bit'' idea is used}\\[1ex]
            1,&\hbox{\hspace*{4mm}otherwise}\\
            \end{array}
    \right.
\eeq

\noindent
so $u - \log_2 p$ bits are required to code the fraction $(b_1, \dots, b_u)\,$.
One bit is required for the sign, and at least 
$\lceil\log_2 (M-m)\rceil$ %
for the
exponent.  Thus
\beq
u - \log_2 p + 1 + \lceil\log_2 (M-m)\rceil \ \leq \ w \,.
\eeq

\noindent
For a sensible design, equality will hold in (2.6), and $M-m$ will be a 
power of two (or one less if the exponent is coded in one's complement
or a special exponent is reserved for zero, but such minor differences 
are unimportant).  Thus (2.4) gives
\beq
2^{-u} kp \ = \ 2^{1-w}R \,.
\eeq

\noindent
The right side depends only on the word length and the range, so (2.7) 
gives a useful relation between the fraction length $u$ and the base
$\beta \ = \ 2^k\,$.

Many different sustems of the class described here have actually been
used. They include, with various 
word lengths, ranges and rounding (or truncating) rules:

\begin{center}
\begin{tabular}{lll}
                       &$\beta = 2,\; p = 2$ &(e.g., PDP 11-45);\\
                       &$\beta = 2,\; p = 1$ &(e.g., CDC 6400);\\
                       &$\beta = 4$        &(e.g., Illiac II);\\
                       &$\beta = 8$        &(e.g., Burroughs 5500); and\\
                       &$\beta = 16$       &(e.g., IBM 360).
\end{tabular}\\
\end{center}

In some machines, bases other than $\beta$ are used in the arithmetic
unit. For example, in the {\footnotesize ILLIAC III} (Atkins~\cite{2}) multiplication
and division are performed with base~$256$, but numbers are stored with
base~$16$.

\section{Other Systems}

\setcounter{equation}{0}

Morris~\cite{22} suggests using ``tapered'' systems in which the division
of bits between the exponent and the fraction depends on the exponent.
The idea is to have a longer fraction for the (commonly occurring)
numbers with exponents close to zero than for numbers with large exponents. 
We do not consider these interesting systems here.

Brown and Richman~\cite{5} assume that floating-point numbers are 
represented in a computer word with two sign bits and a fixed number
of $q$-state devices, for some fixed $q \geq 2$, and they compare bases 
of the form $q^k\,$.  Although the results of Sections 4 and 5 can be
generalized easily to cover their assumptions, we restrict ourselves
to $q=2$, for this is the only case of practical importance.

Finally, we describe a ``logarithmic'' system that is interesting for
theoretical reasons (see Section~4), although it is impractical
(because of the difficulty of performing floating-point additions).
Let $a$ and $b$ be positive integers which, together with the word
length $w\,$, characterize the system.  The floating-point numbers 
are zero and all nonzero real numbers $x$ such that 
$a \cdot \log_2 \vert x \vert + b$ is one of the integers $1, 2, \dots,
2^{w-1} - 1$. If 
\beq
\lambda(x) \ = \ \left\{ \begin{array}{ll}
              0, \hspace*{1mm} &\hbox{if } x = 0\\[1ex]
              \hbox{sign} (x)(a \cdot \log_2 \vert x \vert + b), 
                      \hspace*{1mm}&\hbox{if } x \neq 0\\
            \end{array}
    \right.
\eeq

\noindent
then the floating-point number $x$ may be represented in a computer
word by a convenient code for the integer $\lambda(x)$.  Since
\beq
\hspace*{-9mm} x(\lambda) \ = \ \left\{ \begin{array}{ll}
              0, \hspace*{1mm} &\hbox{if } \lambda = 0\\[1ex]
              \hbox{sign} (\lambda) \cdot 2^{(\lambda - b)/a}, 
                      \hspace*{1mm} &\hbox{if } \lambda \neq 0\\
            \end{array}
    \right.
\eeq

\noindent
the largest and smallest positive floating-point numbers are
\beq
\hspace*{-2mm}f_{\max} \ = \ 2^{(2^{w-1}-1-b)/a}
\eeq
\noindent
and
\beq
\hspace*{-8mm}f_{\min} \ = \ 2^{(1-b)/a}\,,
\eeq

\noindent
respectively, and the range $\log_2 (f_{\max}/f_{\min})$ is
\beq
\hspace*{-1mm}R \ = \ \frac{2^{w-1} - 2}{a} \,.
\eeq

\noindent
For example, taking $a=2^{w-10}$ and $b=2^{w-2}$ gives
$f_{\max} \simeq 2^{256}$, $f_{\min} \simeq 2^{-256}$, and $r \simeq 512$.

If $x$ and $y$ are positive floating-point numbers with $f_{\min}
\leq xy \leq f_{\max}$ then, from (3.1)
\beq
\lambda(xy) \ = \ \lambda(x) + \lambda(y) - b \,.
\eeq

\noindent
Thus, floating-point multiplication and division are easy to perform
in a logarithmic system, and do not introduce any rounding errors.  
Unfortunately, there does not seem to be any easy way to perform
floating-point addition.

\section{The Worst Case Relative-Error Criterion}

\setcounter{equation}{0}

One measure of the precision of a floating-point number system is the
worst relative error $\epsilon$ made in approximating a real number
$x$ (not too large or small) by ${\fl}(x)$, i.e.,
\beq
\displaystyle
\epsilon \ = \ \sup_{f_{\min}\leq\vert x \vert \leq f_{\max}}\,
           \Biggl\vert \frac{x - {\fl}(x)}{x} \Biggr\vert \,. 
\eeq

\noindent
The ``worst case relative-error'' criterion is simply to choose a
number system (with the prescribed $R$ and $w$) to minimise $\epsilon\,$.

For the logarithmic systems described in Section~3, we see from (3.2)
that
\beq
\epsilon \ = \ 2^{1/(2a)} - 1 \ = \ \frac{\log 2}{2a} \,.
\eeq

\noindent
(Here and later we neglect terms of order $a^{-2}$ or $2^{-2u}$,
and logarithms are natural unless otherwise indicated.)  If
\beq
\epsilon_0 \ = \ R2^{-w}\log 2
\eeq
\noindent
then (3.5) and (4.2) give
\beq
\hspace*{-10mm}\epsilon \ = \ \epsilon_0\,.
\eeq

Now consider any floating-point number system with range $R$ and
word length $w\,$.  If $\epsilon$ is defined by (4.1), then
\beq
\hspace*{-10mm}\epsilon \ \geq \ \epsilon_0\,.
\eeq

\noindent
(In a logarithmic system, the logarithms of positive floating-point
numbers are uniformly spaced and all bit patterns are used.)  Thus
we use logarithmic systems as a standard of comparison for other,
more practical, systems.

Wilkinson~\cite{28} shows that
\beq
\epsilon \ = \ 2^{k-u-1}
\eeq

\noindent
for the number systems of Section~2.  From (2.7), (4.3) and (4.6)
\beq
\frac{\epsilon}{\epsilon_0} \ = \ \frac{2^k}{kp \log 2} \ = \ f_1(k,p)
\eeq

\noindent
which shows how much $\epsilon$ exceeds the best possible value
$\epsilon_0$ for a number system with the same $R$ and $w\,$.
Table 1 gives $f_1(k,p)$ for $k=1, 2, \dots, 8\,$.

\begin{center}
\begin{tabular}{cp{2mm}cp{2mm}cp{2mm}cp{2mm}cp{2mm}}
\multicolumn{9}{c}{TABLE 1}\\
\multicolumn{9}{c}{\scriptsize THEORETICAL WORST CASE AND RMS ERRORS*}\\[1ex] \hline\hline
\\[-2ex]
$k$&&$p$&&$\beta=2^k$&&$f_1(k,p)$&&$f_2(k,p)$ \\[1ex] \hline
1&&2&&\W\W2&& \W1.44&&\W1.06\\
1&&1&&\W\W2&& \W2.89&&\W2.12\\
2&&1&&\W\W4&& \W2.89&&\W1.68\\
3&&1&&\W\W8&& \W3.85&&\W1.87\\
4&&1&&\W16&& \W5.77&&\W2.45\\
5&&1&&\W32&& \W9.23&&\W3.51\\
6&&1&&\W64&& 15.4\W&&\W5.34\\
7&&1&&128&& 26.4\W&&\W8.47\\
8&&1&&256&& 46.2\W&&13.9\W\\\hline
\end{tabular}\\[1ex]
* \small See (4.7) and (5.8) for definitions of $f_1$ and $f_2\,$.
\end{center}

The table shows that the implicit-first-bit base~$2$ systems are the best of
those described in Section~2, and close to the best possible, on the
worst case criterion.  Of the explicit-first-bit systems, base~$2$ and
base~$4$ are equally good.  This may be explained as follows.  Changing
from base~$2$ to base~$4$ frees a bit from the exponent for the fraction.
If the first 4-ary digit $d_1$ is 2 or 3, the first fraction bit 
$b_1$ is 1, and the extra fraction bit may increase the precision.
However, if $d_1$ is 1 then $b_1$ is 0, and the bit gained is wasted.
(According to Richman~\cite{24}, Goldberg observed this independently.)
If ${\fl}(x)$ is defined by truncation rather than rounding then
$\epsilon$ is doubled, but the comparison between different bases 
is not changed.

\section{The RMS Relative-Error Criterion}

\setcounter{equation}{0}

Consider forming the product of nonzero floating-point numbers
$x_0, \dots, x_n$ (in one of the usual systems) by $n$ floating-point
multiplications, i.e., define $p_0 = x_0$ and $p_i = {\fl}(p_{i-1}x_i)$
for $i=1, \dots, n\,$. If $\delta_i = (p_{i-1}x_i-p_i)/(p_{i-1}x_i)$
is the relative error made in forming the $i$th product, then the 
relative error in the final result is
\begin{eqnarray}
\displaystyle
\Delta \ &=& \ \frac{x_0 \cdots x_n - p_n}{x_0 \cdots x_n} \ = \
        1 - \prod^{n}_{i=1}\,(1-\delta_i) \nonumber \\
         &=& \ \sum^{n}_{i=1} \delta_i + \hbox{higher order terms}\,.
\end{eqnarray}
\noindent Thus
\beq
\vert \Delta \vert \leq n\epsilon
\eeq

\noindent
where $\epsilon$ is defined by (4.1), and we have neglected a term of
order $n^2\epsilon^2\,$.  Many other bounds on the rounding errors in
algebraic processes are also of the form $f(n)\epsilon$ (see 
Wilkinson~\cite{28},~\cite{29}), which is a good reason for choosing
a floating-point number system according to the worst case criterion
of Section~4.  However, the bound (5.2) is rather pessimistic, for
the individual rounding errors $\delta_i$ in (5.1) usually tend to 
cancel rather than to reinforce each other.  (We are assuming an
unbiased rounding rule as described in Section~6.  With truncation
or biased rounding the bound (5.2) may be realistic.)

If the $\delta_i$ were independent random variables, distributed
with mean 0 and variance $\sigma^2_i\,$, then $\Delta$ would be
distributed with mean 0 and variance $\Sigma^n_{i=1} \sigma^2_i\,$.
Thus a reasonable probabilistic measure of the precision of a 
floating-point number system is the root-mean-square (rms) value
$\delta_{\rm rms}$ of $\delta = (x-{\fl}(x))/x\,$, where $x$ is
distributed like the nonzero results of arithmetic operations performed 
during a typical floating-point computation.

The simulations described in Sections 6 and 7 suggest that the rms
rounding error in floating-point comuptations involving many arithmetic
operations is often roughly proportional to $\delta_{\rm rms}$
(see also \hbox{Weinstein}~\cite{27}).  Thus we prefer $\delta_{\rm rms}$
to other probabilistic measures of precision such as the expected
value of $\vert \delta \vert$ (McKeeman~\cite{21}), the expected
value of $\log_2 \vert \delta \vert$ (Kuki and Codi~\cite{18}), 
and the expected error in ``units in the last place'' (Kahan~\cite{15}).
We disregard errors in the conversion from internal floating-point 
results to decimal output, for the rms value of these errors depends
on the number of decimal places rather than on the internal number
system.  (For the effect of repeated conversions back and forth, 
see Matula~\cite{20}.)

What distribution should we assume for the nonzero real numbers $x$
that are to be approximated by floating-point numbers?  Hamming~\cite{11},
Knuth~\cite{17}, and others argue that we  should assume that
$\log \vert x \vert$ is uniformly distributed.  There are two reasons
why this assumption is only an approximation.  Although $\log \vert x \vert$
may be approximately uniform locally, it is certainly not uniform on
the entire interval $[\log f_{\min}, \log f_{\max}]\,$.  Also, the fine
structure of the distribution is not uniform, for the numbers arising
from multiplications or (more importantly) additions of floating-point
numbers are really discrete rather than continuous variables.  Nevertheless,
we shall make Hamming's assumption in this section.  It is certainly a
much better approximation than assuming that $x$ is uniformly distributed
on some interval.

For the logarithmic systems, $\delta$ is uniformly distributed on
$[-\epsilon_0,\epsilon_0]\,$, where $\epsilon_0$ is given by (4.3).
Thus $\delta_{\rm rms} = \delta_0\,$, where
\beq
\delta_0 \ = \ \frac{\epsilon_0}{\sqrt{3}} \ 
                     = \ \frac{R \cdot \log 2}{2^w \sqrt{3}} \,.
\eeq

\noindent
Because the assumption that $\log \vert x \vert$ is uniform is only
an approximation, there is no result corresponding to the inequality
(4.5), but the logarithmic systems still provide a convenient standard
of comparison for other, more practical, systems.

For the systems of Section~2, there is no loss of generality in assuming
that $x$ lies in $[1/\beta, 1)$ and (by our assumption) $\log_{\beta}x$
is uniformly distributed on $[-1,0)\,$.  Consider numbers $y$
distributed uniformly on a small interval near $x$.  The absolute
error $y - {\fl}(y)$ is approximately uniform on $(-2^{-u-1}, 2^{-u-1})\,$.
(It is certainly not logarithmically distributed, as is assumed to derive 
$(18^{\prime})$ in Benschop and Ratz~\cite{3}.)  Hence
$\alpha = (y - {\fl}(y))/y$ is uniform on $(-2^{-u-1}/x, 2^{-u-1}/x)\,$,
and has probability density function (Feller~\cite{9})
\beq
g_x(\alpha) \ = \ \left\{ \begin{array}{ll}
      2^u x,& \hbox{\hspace*{4mm}if } \vert \alpha \vert < 2^{-u-1}/x\\[1ex]
       0,  & \hbox{\hspace*{4mm}otherwise. }\\
                      \end{array}
    \right.
\eeq

\noindent 
Integrating over the interval $[1/\beta, 1)\,$, we see that $\delta$ is
distributed with density
\begin{eqnarray}
\displaystyle
f(\delta) \ &=& \ \int^1_{x=1/\beta} g_x{(\delta)}\; d\log_{\beta}x \\[2ex]
             &=& \ \left\{ \begin{array}{ll} 
                 2^u(1-2^{-k})/(k \cdot \log 2), \hspace*{10mm}&\hbox{if } \vert \delta \vert < 2^{-u-1}\\[2ex]
                 \left({\displaystyle \frac{1}{2 \vert \delta \vert}}
			- 2^{u-k}\right)\Big/ (k \cdot \log2),
                         \hspace*{6mm}&\hbox{if } 2^{-u-1} \leq \vert \delta \vert < 2^{k-u-1}\\[2ex]
                  0, \hspace*{42mm}&\hbox{otherwise}\,.
                         \end{array}
                  \right.
\end{eqnarray}

\noindent
It is easy to find the expected value of $\delta$, $\delta^2$, 
$\vert \delta \vert$,
$\log_2 \vert \delta \vert$, etc.~from (5.6).  In particular, we find that
$\delta$ is distributed with standard deviation
\beq
\displaystyle
\hspace*{-20mm}\delta_{\rm rms} \ = \ 2^{-u} \sqrt{\frac{4^k-1}{24k \cdot \log 2}}
\eeq

\noindent
and mean 0.  (The mean is actually of order $2^{-2u}\,$, but terms of
this order have been neglected.)

{From} (2.7), (5.3) and (5.6),
\beq
\displaystyle
\frac{\delta_{\rm rms}}{\delta_0} \ = \ \sqrt{\frac{4^k-1}{2p^2(k \cdot \log 2)^3}}
                         \ = \ f_2(k,p)\,,
\eeq

\noindent
and the last column of Table 1 gives $f_2(k,p)$ for $k = 1, \dots, 8\,$.
The table shows that the implicit-first-bit base~$2$ systems are the best of
the systems of Section~2 (and only 6 per cent worse than the logarithmic 
systems) on the rms relative-error criterion.  Base~$4$ (closely followed 
by base~8) is best in the explicit-first-bit systems.  The reason why
base~$4$ is better than explicit base~$2$ is apparent from the discussion
at the end of Section~4:  $\vert \delta \vert$ is never greater for 
base~$4$ than for explicit base~$2$, and sometimes it is smaller.  A similar
argument shows that implicit base~$2$ is better than base~$4$.

Because of the different ranges possible with base~$4$ and base~$8$, there 
are some choices of minimal acceptable range for which base~$8$ is 
preferable to base~$4$, but bases higher than 8 are always inferior to
base~$4$ on the rms relative-error criterion.

\section{Simulation of Different Systems}

\setcounter{equation}{0}

Three classes of floating-point computations were run, using various number
systems with $w = 32$ and $R \simeq 512$ (the same as for single-precision
on the IBM 360 and many other computers).  The systems were a logarithmic
system $S_0$ with $a = 2^{22}$ and $b = 2^{30}$ (see Section~3), and the
following examples of the systems described in Section~2.\\

\begin{tabular}{lll}
$S_1:$&$ \beta=2,\; u=23,\; p=2$ &(base~$2$ with a 23-bit fraction, the first 
                             bit implicit).\\
$S_2:$&$ \beta=4,\; u=23$ &(base~$4$ with 23 bits or $11\frac{1}{2}$ digits).\\
$S_3:$&$ \beta=2,\; u=22,\; p=1$ & (base~$2$ with 22 bits, all explicit).\\
$S_4:$&$ \beta=16,\; u=24$ &(base~$16$ with 24 bits or 6 digits).\\
$S^{'}_4: $&\multicolumn{2}{l}{The same as $S_4$ with truncation (towards zero) rather than
                         rounding.}\\
$S_5:$&$ \beta= 256,\; u=25$ &(base~$256$ with 25 bits or $3\frac{1}{8}$ digits).
\end{tabular}\\

The rounding rule for systems $S_1$ to $S_5$ is the ``$R$*-mode'' of
Kuki and Cody~\cite{18}: ${\fl}(x)$ is defined to be the floating-point
number closest to $x$, and ties are broken by choosing ${\fl}(x)$ so that
its least significant fraction bit is one.  Formally, if $x$ is a 
nonzero real number with binary expansion
\beq
x = s \sum^\infty_{j=1} b_j 2^{ke-j}
\eeq

\noindent
(taking the terminating expansion if there is one, normalizing so that
one of $b_1, \dots, b_k$ is nonzero, and neglecting the possibility
of underflow or overflow), then 

\beq
{\fl}(x) \ = \ \left\{ \begin{array}{ll} \displaystyle
s \sum^u_{j=1} b_j 2^{ke-j}, \hspace*{20mm} &\hbox{if \ } b_{u+1} = 0 
\hbox{ \ or }\sum^\infty_{j=0} b_{u+j}2^{-j}=\frac{3}{2} \nonumber\\[4ex]
\displaystyle s \Bigg(\sum^u_{j=1} b_j 2^{ke-j}+2^{ke-u}\Bigg)\,, 
                        \hspace*{2mm}&\hbox{otherwise}\,. \nonumber
            \end{array}
    \right.
\eeq

\vspace*{3mm}
\noindent
The special case $b_u b_{u+1} \cdots = 11000 \cdots$ is quite important,
for it often occurs when $x$ is the result of a floating-point addition,
and neglecting it can lead to bias in the rounding.

All the floating point number systems were simulated on an IBM 360/91
computer, with arithmetic operations performed in double precision
$(\beta = 16, u = 56)$ before rounding or truncating approximately.
Thus, the number of guard units used was effectively infinite.  The 
data were pseudorandom double-precision numbers distributed as described 
in Section~7, and ``exact'' results were computed using double precision
throughout.

Forming sums, solving systems of linear equations, and finding the
eigenvalues of symmetric matrices were the chosen classes of floating-point
computations.  They appear to be fairly typical of computations in which 
the effect of rounding errors may be important.  Details, and the results
of the simulations, are given in Section~7.  Other classes that have been
considered include solving ordinary differential equations 
(Henrici~\cite{12},~\cite{13}, Hull and Swenson~\cite{14}),
fast Fourier transforms (Kaneko and Liu~\cite{16}, Ramos~\cite{23},
and Weinstein~\cite{27}), matrix iterative processes (Benschop and
Ratz~\cite{3}), solving positive-definite linear systems (Tienari~\cite{25}),
and forming products (Section~5).

\section{Details and Results of the Simulations}

\setcounter{equation}{0}

\noindent{\large\it Sums}\\[-1ex]

Let $m$ and $n$ be positive integers.  A number $z$ was drawn from a
uniform distribution on $[0,1]$, then numbers $x_1, \dots, x_n$
were drawn independently from a uniform distribution on $[-Z, Z]$,
where $Z = 256^z$ is a scale factor used to avoid a bias in favour
of any of the number systems (see Kuki and Cody~\cite{18}).  The
approximate sums $s_j$ of ${\fl}(x_1), \dots, {\fl}(x_n)$ were accumulated,
in the usual way, with each of the number systems $S_j$ described in
Section~6, and the errors
\beq
\alpha_j \ = \ \frac{\displaystyle \sum^n_{i=1}x_i\;\; - \;\;s_j}
                    {\displaystyle \sum^n_{i=1} \vert x_i \vert}
\eeq

\noindent
were found.  (The denominator is used in preference to
$\sum^n_{i=1}x_i$ to ensure that $\alpha_j$ is small.)  The procedure
was repeated $m$ times and the rms values $\beta_j$ of the $\alpha_j$
were found.  For purposes of comparison between the systems, it is 
convenient to consider the normalized rms errors 
$\gamma_j=\beta_j/\beta_0\,$.  (Recall that $\beta_0$ is the rms error
for the logarithmic system $S_0\,$.)

Table 2 gives $\gamma_j$ for various choices of $m$ and $n\,$.  If the
$\alpha_j$ are considered as random variables drawn from a distribution
with mean square $B^2_j\,$, then $\beta_j$ and $\gamma_j$ may be
regarded as estimates of $B_j$ and $B_j/B_0\,$, respectively.  $m$
was chosen large enough to ensure that the standard error of the
estimates $\gamma_j$ given in Tables 2--4 is less than five units
in the last decimal place.

For $n=1$ we are merely estimating the rms relative error in
approximating $x_1$ by ${\fl}(x_1)\,$, and the results agree with the
predictions of Section~5 (see the last column of Table 1).
Except for $S'_4\,$, the effect of varying $n$ is small, and does
not affect the ranking of the systems.

It may be shown that
\beq
B_j \ = \ \left\{ \begin{array}{ll}
                 O(n^{3/2})\,, & \hspace*{6mm}\hbox{for } S'_4 \\[1ex]
                 O(n)\,,       & \hspace*{6mm}\hbox{for the other systems}\\
               \end{array} \right.
\eeq

\noindent
so it is not surprising that $\gamma'_4$ appears
to grow like $n^{1/2}\,$.  (The same applies if truncation is
downwards instead of towards zero.)  Results for sums of positive 
numbers are similar, although $B_j$ is larger by a factor of order
$n^{1/2}$ for all the systems.

\begin{center}
\begin{tabular}{cp{1mm}cp{1mm}cp{1mm}cp{1mm}cp{1mm}cp{1mm}cp{1mm}c}
\multicolumn{15}{c}{TABLE 2}\\
\multicolumn{15}{c}{\scriptsize RESULTS FOR SUMS}\\[1ex] \hline\hline
&& && && && && && &&\\[-2ex]
$n$&&$m/1000$&&$\gamma_1$&&$\gamma_2$&&$\gamma_3$&&$\gamma_4$&&$\gamma'_4$&&
                                      $\gamma_5$ \\[1ex] \hline
&& && && && && && &&\\[-2ex]
\W\W1 &&1000 &&1.06     &&1.68 &&2.12 &&2.45 &&\W4.89 &&13.9 \\
\W\W2 &&\W100 &&1.11    &&1.68 &&2.23 &&2.38 &&\W5.53 &&13.4 \\
\W\W4 &&\W100 &&1.13    &&1.69 &&2.25 &&2.36 &&\W6.33 &&13.2 \\
\W\W8 &&\W100 &&1.12    &&1.69 &&2.24 &&2.36 &&\W7.95 &&13.2 \\
\W10 &&\W100 &&1.12   &&1.69 &&2.23 &&2.36 &&\W8.76 &&13.4 \\
\W16 &&\W\W10 &&1.11  &&1.72 &&2.22 &&2.37 &&10.9\W &&13.3 \\
\W32 &&\W\W10 &&1.09 &&1.71 &&2.18 &&2.39 &&15.9\W &&13.6 \\
\W64 &&\W\W10 &&1.08 &&1.67 &&2.14 &&2.43 &&22.4\W &&13.9 \\
100 &&\W\W30 &&1.06 &&1.68 &&2.13 &&2.41 &&28.1\W &&13.6 \\[1ex] \hline

\end{tabular}
\end{center}

\vspace*{3mm}

\noindent{\large\it Solving Systems of Linear Equations}\\[-1ex]

$z_1$ and $z_2$ were drawn independently from a uniform distribution
on $[0, 1]\,$, giving scale factors $Z_1 = 256^{z_1}$ and 
$Z_2 = 256^{z_2}\,$.  
Numbers $a_{p,q}\;(p,q=1,\dots,n)$ were drawn independently
from a uniform distribution on $[-Z_1, Z_1]\,$; and $x_1,\dots,x_n$
were drawn similarly from $[-Z_2, Z_2]\,$.  For each of the number
systems $S_j\,$, let $A^{(j)}=({\fl}(a_{p,q}))$, $A=(a_{p,q})$, 
$x=(x_p)$, $b = (b_p)=Ax$,
and $b^{(j)}=({\fl}(b_p))$.  The system of equations
\beq
A^{(j)}y \ = \ b^{(j)}
\eeq

\noindent
was solved by Gaussian elimination with complete pivoting, giving the
approximate solution $y^{(j)}\,$, and the error
\beq
\alpha_j \ = \ \frac{\Vert Ay^{(j)} - b \Vert_2}
                    {\Vert A \Vert_E \ \Vert x \Vert_2}
\eeq

\noindent
was computed.  (Here $\Vert A \Vert_E = \Big(\sum^n_{p=1} \sum^n_{q=1}
a^2_{p,q}\Big)^{1/2}\,$.  From results of Wilkinson~\cite{28},
\cite{29}, $\alpha_j$ is small even if $A$ is rather ill conditioned.)
The procedure was repeated $m$ times, the rms values $\beta_j$
of the $\alpha_j$ were computed, and the ratios $\gamma_j
=\beta_j/\beta_0$ were found.  The results for various $m$ and $n$
are given in Table 3. 

\begin{center}
\begin{tabular}{cp{1mm}cp{1mm}cp{1mm}cp{1mm}cp{1mm}cp{1mm}cp{1mm}c}
\multicolumn{15}{c}{TABLE 3}\\
\multicolumn{15}{c}{\scriptsize RESULTS FOR SYSTEMS OF LINEAR
                                    EQUATIONS}\\[1ex] \hline\hline
&& && && && && && &&\\[-2ex]
$n$&&$m/1000$&&$\gamma_1$&&$\gamma_2$&&$\gamma_3$&&$\gamma_4$&&$\gamma'_4$&&
                                      $\gamma_5$ \\[1ex] \hline
&& && && && && && &&\\[-2ex]
\W\W1 &&100 &&1.30         &&2.06 &&2.61 &&2.99 &&4.92 &&17.0 \\
\W\W2 &&100 &&1.30         &&2.01 &&2.59 &&2.90 &&5.33 &&16.3 \\
\W\W4 &&\W10 &&1.27        &&1.97 &&2.56 &&2.80 &&5.63 &&15.7 \\
\W\W8 &&\W\W4 &&1.23       &&1.89 &&2.45 &&2.65 &&6.1\W &&14.9 \\
\W16 &&\W\W1 &&1.18        &&1.82 &&2.35 &&2.60 &&7.1\W &&14.4\\[1ex] \hline

\end{tabular}
\end{center}

Multiplication and division are performed exactly in a logarithmic
system, so $\beta_0$ is less than would otherwise be expected, and
$\gamma_1,\dots,\gamma_5$ are higher than for sums, especially for
small values of $n\,$.  However, the ratios of $\gamma_1,\dots,\gamma_5$
are much the same as for sums, and the ranking of the systems is
preserved.  Results for positive $a_{p,q}$ and/or $x_p$ are similar.

\pagebreak[3]
{\samepage
It is interesting that $\gamma_4 < 2\gamma'_4$ for $n=1$ and 2.
When $n=1$ and $S'_4$ is used, the errors made in forming ${\fl}(a_{1,1})$
and ${\fl}(b_1)$ tend to cancel when ${\fl}(b_1)/{\fl}(a_{1,1})$ is computed.
Presumably there is a similar, though less marked, effect for $n>1\,$.\\
}

\noindent
{\large\it Finding Eigenvalues of Symmetric Matrices}\\

Numbers $a_{p,q}\; (1 \leq p \leq q \leq n)$ were drawn independently
from a uniform distribution on $[-Z,Z]\,$, where $Z$ was a scale factor
chosen as above.  The other elements of $A = (a_{p,q})$ were defined 
by symmetry.  For each number system $S_j\,$, the approximate
eigenvalues $\lambda^{(j)}_1 \leq \cdots \leq \lambda^{(j)}_n$ of
$A^{(j)} = ({\fl}(a_{p,q}))$ were computed by reducing $A^{(j)}$ to
tridiagonal form and then using the QR algorithm (Wilkinson~\cite{29}).
We used translations of the Algol 60 procedures TRED1 
(Martin {\etal}~\cite{19}) and TQL1 (Bowler {\etal}~\cite{4}), except
for some trivial modifications to avoid unnecessary rounding errors
when $n=2\,$.  The stopping criterion for the QR algorithm was the same 
for all number systems.  (The parameters {\it macheps\/} and {\it tol\/}
of the procedures were set to $10^{-8}$ and $10^{-60}\,$, respectively.)
The errors
\beq
\displaystyle
\alpha_j \ = \ \Bigg( \sum^n_{i=1} (\lambda_i - \lambda^{(j)})^2\Bigg)^{\frac{1}{2}}\Bigg/ \
                    \Vert A \Vert_E
\eeq 

\noindent
were computed.  (Here $\lambda_1 \leq \cdots \leq \lambda_n$ are the
exact eigenvalues of $A\,$.)  The procedure was repeated $m$ times,
the rms values $\beta_j$ of the $\alpha_j$ computed, and the ratios
$\gamma_j=\beta_j/\beta_0$ found, as above.  The results are given
in Table 4. 

\begin{center}
\begin{tabular}{cp{1mm}cp{1mm}cp{1mm}cp{1mm}cp{1mm}cp{1mm}cp{1mm}c}
\multicolumn{15}{c}{TABLE 4}\\
\multicolumn{15}{c}{\scriptsize RESULTS FOR EIGENVALUES OF SYMMETRIC
                                    MATRICES}\\[1ex] \hline\hline
&& && && && && && &&\\[-2ex]
$n$&&$m/1000$&&$\gamma_1$&&$\gamma_2$&&$\gamma_3$&&$\gamma_4$&&$\gamma'_4$&&
                                      $\gamma_5$ \\[1ex] \hline
&& && && && && && &&\\[-2ex]
\W\W2 &&100   &&1.07 &&1.61 &&2.14 &&2.38 &&\W6.06 &&15.2 \\
\W\W4 &&\W10  &&1.33 &&2.24 &&2.65 &&3.60 &&10.5\W &&25.8 \\
\W\W8 &&\W\W3 &&1.14 &&2.01 &&2.34 &&3.73 &&10.8\W &&29.6 \\
\W16  &&\W\W1 &&1.00 &&1.82 &&1.99 &&3.49 &&10.7\W &&28.8\\[1ex] \hline

\end{tabular}
\end{center}

The method used for finding eigenvalues depends heavily on multiplications
by matrices of the form {\footnotesize$\pmatrix{c & s \cr -s & c}\,$}, where $c^2 +
s^2 = 1\,$. The numbers $c$ and $s$ are certainly not distributed as
assumed in Section~5.  This, along with other observations made above,
may explain the interesting variations in the $\gamma_j\,$.  Despite
these variations, the ranking of the different systems is as predicted in
Section~5.

\section{Conclusions}

Comparing $\gamma'_4$ with $\gamma_4$ in Tables 2--4 shows that the rms
error for truncation is usually considerably more than twice as much
as for rounding.  However, truncation is often preferred because the
usual implementation of rounding requires an extra carry propagation.
An interesting compromise is the ``von Neumann round'' (Burks {\etal}~\cite{6}, 
Urabe~\cite{26}), for which the result of an arithmetic 
operation is truncated, and then the least significant bit is set to one.  
(An exception could be made if the result is exactly representable;
this would involve checking if the truncated bits were all zero.)
No extra carry propagation is required, and the rms error is twice
that for normal rounding, so considerably better than for truncation.

The most accurate practical systems are base~$2$ with the first fraction 
bit implicit.  If the accuracy gained by having the first bit implicit is
not considered sufficient compensation for the disadvantages entailed, 
then base~$4$ (or perhaps base~$8$) is the best choice.

The accuracy lost by using base~$16$ or higher is roughly as predicted
in Section~5.  High bases may have some implementation advantages
(Anderson {\etal}~\cite{1}, Atkins~\cite{2}).  In practice both factors
should be considered.  The number of guard digits used is also important.
The use of high bases, only one guard digit, and truncation instead
of rounding is probably acceptable on machines with a long floating-point
word.  However, to minimize the need for double-precision computations,
it seems wise to try to squeeze out the last drop of accuracy on a
computer with a short floating-point word (say 32--40 bits).  The
amount that can be squeezed out is often significant.  For example,
our simulations show that using system $S_1$ instead of $S'_4$ is
roughly equivalent to carrying one more {\it decimal\/} place.

\section*{Acknowledgement} 

The author wishes to thank W J Cody, the late Prof.~G E Forsythe,
and Prof.~W Kahan for their comments on an earlier version of this
paper.

\end{document}